\documentclass[11pt]{article}
\usepackage{amsmath}
\usepackage{amscd}
\usepackage{amsfonts}
\usepackage{amssymb}
\usepackage{amsthm}

\newcommand{\rar}{\rightarrow}
\newcommand{\calg}{\mathcal}

\newcommand{\rhh}{{\text{R}}\text{Hom}}
\newcommand{\rhhh}{\text{R}\mathcal{H}om}

\newcommand{\strc}{\mathcal{O}_{X}}
\newcommand{\strcc}[1]{\mathcal{O}_{#1}}

\newcommand{\dcat}[1]{\text{D}^{b}(#1)}

\newcommand{\hkr}{\text{I}_{HKR}}

\newcommand{\hoch}[2]{\text{HH}_{#1}(#2)}
\newcommand{\hochh}[1]{\text{HH}_{\bullet}(#1)}

\newcommand{\per}[1]{\text{perf }({#1})}
\newcommand{\natl}[1]{{#1}^{\text{nat}}_*}
\newcommand{\muk}[1]{{#1}^{\text{muk}}_*}
\newcommand{\cald}[1]{{#1}^{\text{cal}}_*}
\newcommand{\gpair}{\langle \text{ },\text{ }\rangle}
\newcommand{\spair}[2]{\langle #1,#2 \rangle}
\newcommand{\splt}{\Phi_{\text{pt} \rightarrow X \times Y}}
\newcommand{\perr}[1]{\text{Z}^0(\text{perf}({#1}))}

\newtheorem{thm}{Theorem}

\newtheorem{prop}{Proposition}

\title{The Mukai pairing and integral transforms in Hochschild homology.}

\author{Ajay C. Ramadoss}

\begin{document}

\maketitle

\begin{abstract}
Let $X$ be a smooth proper scheme over a field of characteristic
$0$. Following Shklyarov [10] , we construct a (non-degenerate)
pairing on the Hochschild homology of $\per{X}$, and hence, on the
Hochschild homology  of $X$. On the other hand the Hochschild
homology of $X$ also has the Mukai pairing (see [1]). If $X$ is
Calabi-Yau, this pairing arises from the action of the class of a
genus $0$ Riemann-surface with two incoming closed boundaries and no
outgoing boundary in $\text{H}_{0}({\mathcal M}_0(2,0))$ on the
algebra of closed states of a version of the B-Model on $X$. We show
that these pairings "almost" coincide. This is done via a different
view of the construction of integral transforms in Hochschild
homology that originally appeared in Caldararu's work [1]. This is
used to prove that the more "natural" construction of integral
transforms in Hochschild homology by Shklyarov [10] coincides with
that of Caldararu [1]. These results give rise to a Hirzebruch
Riemann-Roch theorem for
the sheafification of the Dennis trace map. \\
\end{abstract}

\section*{Introduction.}

Let $X$ be a smooth proper scheme over a field $\mathbb K$ of
characteristic $0$. Let $\per{X}$ denote the DG-category of left
bounded perfect injective complexes of $\strc$-modules . There is a
natural isomorphism of Hochschild homologies (see [5] for instance)
\begin{equation} \label{one} \text{HH}_{\bullet}(X) \simeq \text{HH}_{\bullet}(\per{X}) \text{ . }\end{equation}

If $Y$ is any smooth proper scheme, an object $\Phi \in \per{X
\times Y}$ can be thought of as the kernel of an integral transform
from $\per{X}$ to $\per{Y}$ (Section 8 of [11]). This is a morphism
from $\per{X}$ to $\per{Y}$ in the homotopy category
$\text{Ho}(\text{dg-cat})$ of dg-categories modulo
quasi-equivalences. We will abuse notation and denote this by $\Phi$
as well. It follows that $\Phi$ induces a map
$\Phi_*:\text{HH}_{\bullet}(\per{X}) \rar
\text{HH}_{\bullet}(\per{Y})$ and hence, by $\eqref{one}$ , a map
$$\Phi^{\text{nat}}_*: \text{HH}_{\bullet}({X}) \rar
\text{HH}_{\bullet}({Y}) \text{ . }$$

One also has (see [10] ) a Kunneth quasiisomorphism
$$K:\hochh{\per{X}} \otimes \hochh{\per{Y}} \rar \hochh{\per{X
\times Y}} \text{ . }$$ Since $X$ is smooth, the diagonal $\Delta: X
\rar X \times X$ is a local complete intersection. Hence,
$\strcc{\Delta}:= {\text{R}}\Delta_* \strc$ is a perfect complex on
$X \times X$ (see [11], Section 8). We will abuse notation and
denote $\strcc{\Delta}$ thought of as the kernel of an integral
transform from $X \times X$ to $\text{Spec }\mathbb K$ by $\Delta$.
One then has a pairing given by the composite map
$$ \begin{CD}
\hochh{\per{X}} \otimes \hochh{\per{X}}\\
 @VKVV\\
  \hochh{\per{X \times X}} @>\Delta_*>> \hochh{\per{\mathbb K}}=\mathbb K \text{ .
} \\
\end{CD}$$ Denote this
pairing by $\langle \text{  },\text{  } \rangle_{\text{Shk}}$.\\

On the other hand, the work of A. Caldararu [1] constructs the
following:\\

$\bullet$ A non-degenerate Mukai pairing
$$\langle \text{ },\text{ } \rangle_M:\hochh{X} \otimes \hochh{X}
\rar \mathbb K
\text{ . }$$  \\
$\bullet$ For each $\Phi \in \text{Perf }(X \times Y)$ an "integral
transform" $$\Phi^{\text{cal}}_*:\hochh{X} \rar \hochh{Y} \text{ .
}$$

If $X$ is Calabi-Yau, it has been argued implicitly by Caldararu [3]
that $\langle \text{ },\text{  }\rangle_M$ is precisely the pairing
on
 $\hochh{X}$ arising from the action (on $\hochh{X}$) of the class of a
genus $0$ Riemann-surface with two incoming closed boundaries and no
outgoing boundary in $\text{H}_{0}({\mathcal M}_0(2,0))$. Let
$\vee:\hochh{X} \rar \hochh{X}$ be the whose image under the
Hochschild-Kostant-Rosenberg isomorphism is the involution on Hodge
cohomology that acts on the direct summand
$\text{H}^{q}(X,\Omega_X^p)$ by multiplication by${(-1)}^p$.

\subsection*{The "natural" pairing and the Mukai pairing.}

The main result of this note is as follows.

\begin{thm}  Let $a,b \in \hochh{X}$. Then,
$$\langle b^{\vee},a \rangle_M = \langle a, b \rangle_{\text{Shk}}
\text{ . }$$
\end{thm}
If $X$ is a smooth proper quasi-compact scheme, the category
$\per{X}$ is quasi-equivalent to $\per{A}$ for some DG-algebra $A$
(see [6],[11]). In this case, the pairing $\langle \text{ },\text{
}\rangle_{\text{Shk}}$ on $\hochh{X}$ is the pairing on $\hochh{A}$
described in [10]. On the other hand, the Mukai pairing
$\langle\text{ },\text{  }\rangle_M$ has been explicitly computed at
the level of Hodge cohomology in [8]. In an implicit form, this
computation appeared earlier in [7]. Theorem 1 therefore, enables us
to relate the familiar Riemann-Roch-Hirzebruch theorem for a proper
scheme over $\mathbb K$ to the more abstract "noncommutative"
Riemann-Roch theorem in [10].\\

Further, if $X$ is Calabi-Yau, so is $A$. In this case Theorem 1 is
very similar to Conjecture 6.2 in [10] for proper homologically
smooth Calabi-Yau DG-algebras $A$ such that $\per{A}$ is
quasi-equivalent to $\per{X}$ for some smooth proper quasi-compact
scheme $X$. We make a remark about
this in Section 2.3.\\

\subsection*{Integral transforms in Hochschild homology.}

Let us outline how Theorem 1 is proven. It was stated and proven in
[10] that if $\Phi \in \per{X \times Y}$, then $\natl{\Phi}$ is
simply convolution with the Chern character of $\Phi$ with respect
to the pairing $\gpair_{\text{Shk}}$. Besides [10], the reader may
refer to Theorems 4 and 5 in this paper for the precise statement.
We construct a map $\muk{\Phi}:\hochh{X} \rar \hochh{Y}$ that is
"almost" convolution with the Chern character of $\Phi$ with respect
to the Mukai pairing. We then proceed to prove that $\muk{\Phi}$ has
all the "good properties" one expects of an integral transform in
Hochschild homology (Propositions 1 ,2 and 3 of this paper). We
recall that the integral transform from $\per{X}$ to $\per{X}$
arising out of the element $\strcc{\Delta}$ of $\per{X \times X}$ is
the identity. It follows that $\natl{\strcc{\Delta}}
 =\text{id}$. Proposition 2, which says that $\muk{\strcc{\Delta}} = \text{id}$ as well, is
then used to prove Theorem 1.\\

The fact that $\muk{\Phi}$ has all the "good properties" one expects
of an integral transform in Hochschild homology  is also exploited
to prove the following theorem.

\begin{thm}
$$\natl{\Phi} = \muk{\Phi} = \cald{\Phi} \text{ . }$$
\end{thm}

In other words,the "good constructions" of integral transforms in
Hochschild homology coincide.\\

\subsection*{A Hirzebruch-Riemann-Roch for the sheafification of the
Dennis trace map.}

We now mention another consequence of Theorems 1 and 2. Recall that
we have an isomorphism of higher K groups
$$K_i(X) \simeq K_i(\per{X}) \text{ . }$$
For any DG-category $\calg C$, let $\text{Z}^0(\calg C)$ denote the
category such that $$\text{Obj}(\text{Z}^0(\calg C)) =
\text{Obj}(\calg C) \text{ and } \text{Hom}_{\text{Z}^0(\calg
C)}(M,N) = \text{Z}^0(\text{Hom}_{\calg C}(M,N)) \text{ }$$ $$
\forall \text{ } M,N \in \text{Obj}(\calg C) \text{ . }$$ Here,
$\text{Z}^0(\text{C})$ is the space of $0$-cocycles for any cochain
complex $\text{C}$. If $\text{Z}^{0}(\calg C)$ is exact, one has a
Dennis trace map
$$\text{Ch}^i:K_i(\calg C) \rar \hoch{i}{\calg C}$$
(see [12]). This therefore, yields us a map $$\text{Ch}^i:K_i(X)
\rar \hoch{i}{\per{X}} \simeq \hoch{i}{X} \text{ . }$$ This map is
the "sheafification of the Dennis trace map" constructed in [13].
Let $\hkr: \hochh{X} \rar \oplus_{p,q} \text{H}^p(X,\Omega_X^q)$
denote the Hochschild-Kostant-Rosenberg isomorphism. Let
$$\text{ch}^i:K_i(X) \rar \oplus_j \text{H}^{j-i}(X,\Omega_X^j) $$
denote $\hkr \circ \text{Ch}^i$. It was proven in [2] (Theorem 4.5)
that $\text{ch}^0$ is the usual Chern character. We have the
following generalization of the Hirzebruch Riemann-Roch theorem.

\begin{thm}
Let $f:X \rar Y$ be a smooth proper morphism between proper schemes
$X$ and $Y$. Let $Z$ be a smooth quasi-compact separated scheme.
Then,
$${(f \times \text{id})}_*(\text{ch}^i(\alpha) \pi_X^* \text{td}(T_X)) =
\text{ch}^i({(f \times \text{id})}_*(\alpha)) \pi_Y^* \text{td}(T_Y)$$ for any $\alpha \in K_i(X \times Z)$.\\
\end{thm}

\subsection*{Layout of this note.}

Section 1 reviews some basic facts from D. Shklyarov's work [10].
Section 2.1 recalls A. Caldararu's construction of the Mukai pairing
[1] and related results. In Section 2.2, we give an alternate
construction of $\Phi_*:\hochh{X} \rar \hochh{Y}$ for any $\Phi \in
\per{X \times Y}$. We prove Theorem 1 and Theorem 2 in Section 2.2.
Section 2.3 contains some remarks about what Theorem 1 means when
$X$ is Calabi-Yau. Section 2.4 proves Theorem 3. \\

\subsection*{Acknowledgements.} I am grateful to Prof. Kevin Costello, Prof. Madhav Nori and Prof. Boris Tsygan for some very useful discussions.

\section{The "natural pairing" on the Hochschild homology of schemes.}

This section primarily recalls material from D. Shklyarov's work [10]. The term "DG algebra" in this section shall refer to a proper homologically smooth DG-algebra unless explicitly stated otherwise.\\

\subsection{Preliminary recollections.}
Recall that a DG-algebra $A$ is proper if $\sum_n \text{dim }\text{H}^n(A)  < \infty$
 and is homologically smooth if it is quasi isomorphic to a perfect $A^{\text{op}} \otimes A$-module. Here, $A^{\text{op}}$ denotes the
  opposite algebra of $A$. The term "$A$-module" shall refer to a right $A$-module.\\

 Recall that a $A$-module is said to be semi-free if it is obtained from a finite set of free $A$-modules after taking
 finitely many cones of degree $0$ closed morphisms . A perfect $A$-module is a direct summand of a semi-free $A$-module.
 Let $\per{A}$ denote the DG-category of perfect $A$-modules.
 We recall the following facts from [10].\\

 Fact 1: If $A$ is a DG-algebra, the natural embedding of the category with a unique object whose morphisms are given by $A$
 into $\per{A}$ induces an isomorphism
\begin{equation} \label{onee} \hochh{A} \simeq \hochh{\per{A}} \end{equation}

 Fact 2: If $A$ and $B$ are DG-algebras and $\Phi$ is a perfect $\text{A}^{\text{op}} \otimes B$-module, then $\Phi$ gives a (DG) functor
 $$\Phi_*:\per{A} \rar \per{B}$$
 $$ M \leadsto M \otimes_A \Phi \text{ . }$$
 $\Phi_*$ therefore induces a map
 $$\Phi^{\text{nat}}_*:\hochh{\per{A}} \rar \hochh{\per{B}} \text{ . }$$

 Fact 3: Let $\Delta$ denote $A$ treated as a a perfect $A^{\text{op}} \otimes A$-module in the natural way. Then,
  by Fact 2, we have a DG functor $\Delta_*:\per{A \otimes A^{\text{op}}} \rar \per{\mathbb K}$. Further, there is a isomorphism
 $$K: \hochh{\per{A}} \otimes \hochh{\per{A^{\text{op}}}} \rar \hochh{\per{A \otimes A^{\text{op}}}} \text{ . }$$
 The map $\natl{\Delta} \circ K:\hochh{\per{A}} \otimes \hochh{\per{A^{\text{op}}}} \rar \hochh{\per{\mathbb K}}=\mathbb K$ therefore
 gives rise to a pairing
 $$\langle\text{ },\text{  }\rangle_{\text{Shk}}:\hochh{A} \otimes \hochh{A^{\text{op}}} \rar \mathbb K \text{ . }$$

 For any exact $\mathbb K$-linear category $\mathcal C$, let $K_0({\mathcal C})$ denote the Grothendieck group of $\mathcal C$.
 Recall from [10] that there is a Chern character
 $$\text{Ch}: K_0(\per{A}) \rar \text{HH}_{0}(\per{A}) \simeq \hoch{0}{A} \text{ . }$$
 Let $A$ and $B$ be DG-algebras. We abuse notation and denote the composite map
 $$ \begin{CD} \hochh{A} \otimes \hochh{A^{\text{op}}} \otimes \hochh{B} @> \langle\text{ },\text{ }\rangle_{\text{Shk}} \otimes \text{id} >> \hochh{B} \end{CD} $$
  by $\langle\text{ },\text{ }\rangle_{\text{Shk}}$ itself. Identify $\hochh{A^{\text{op}} \otimes B}$ with $\hochh{A^{\text{op}}} \otimes \hochh{B}$
  via the inverse of the Kunneth isomorphism. If $\Phi \in \per{A^{\text{op}} \otimes B}$,
  the following theorem from [10] (Theorem 3.4 of [10]) says that $\natl{\Phi}$ is just "convolution with $\text{Ch}(\Phi)$".

\begin{thm} $$\natl{\Phi}(x)= \langle x,\text{Ch}(\Phi)\rangle_{\text{Shk}} $$ for any $x \in \hochh{A}$.
\end{thm}

 Note that Theorem 4 implies that $\natl{\Phi}$ depends only on the
 image of $\Phi$ in $\text{D}(\per{A^{\text{op}} \otimes B})$.\\

\subsection{The natural pairing on the Hochschild homology of schemes.}
In this subsection, whenever $f:X \rar Y$ is a morphism of schemes,
$f_*$,$f^*$ etc shall denote the corresponding derived functors. Let
$X$ be a quasicompact separated scheme over $\mathbb K$. In this
case, the (unbounded) derived category $\text{D}_{qcoh}(X)$ of
quasi-coherent $\strc$-modules on $X$ admits at least compact
generator $E$ (see [11]). This is a perfect complex of
$\strc$-modules. We
recall the following facts.\\

{\it Fact 1:} For each compact generator $E$ of $\text{D}_{qcoh}(X)$ there one can choose a (proper if and only if $X$ is proper)
DG-algebra ${\calg A}(E)$ such that $\per{{\calg A}(E)}$ is quasi-equivalent to $\per{X}$
 (see [6],[11 ]). \\

{\it Fact 2:} Recall that if $E$ is a compact generator of $\text{D}_{qcoh}(X)$ and if $F$ is a compact generator
of $\text{D}_{qcoh}(Y)$ then $E \boxtimes F$ is a compact generator of $\text{D}_{qcoh}(X \times_{\mathbb K} Y)$. \\

{\it Fact 3:} The ${\calg A}(E)$ can be chosen so that
$${\calg A}(E \boxtimes F) = {\calg A}(E) \otimes {\calg A}(F)$$ whenever $E$ and $F$ are as in Fact 2 above.\\

{\it Fact 4:} If $E$ is a compact generator of $\text{D}_{qcoh}(X)$, so is the dual perfect complex $E^{\vee}$. One can
choose ${\calg A}(E^{\vee})$ to be ${\calg A}(E)^{\text{op}}$. Hence, $\per{{\calg A}(E)}$ is quasi-equivalent to $\per{{\calg A}(E)^{\text{op}}}$. \\

From the quasi-equivalences $\per{{\calg A}(E)} \simeq \per{X}$ and $\per{{\calg A}(E)^{\text{op}}} \simeq \per{X}$, we obtain isomorphisms
$$i:\hochh{X} \simeq \hochh{{\calg A}(E)}$$ $$j:\hochh{X} \simeq \hochh{{\calg A}(E)^{\text{op}}} \text{ . }$$
For $X$ proper let $\langle \text{ },\text{ } \rangle_{\text{Shk}}$ be the pairing on $\hochh{X}$ such that
$$\langle a,b \rangle_{\text{Shk}} = \langle i(a),j(b) \rangle_{\text{Shk}}$$
for all $a,b \in \hochh{X}$. Note that the RHS of the above equation
has been defined in the previous subsection. We identify $\hochh{X
\times Y}$ with $\hochh{X} \otimes \hochh{Y}$ via the inverse of the
Kunneth isomorphism. Recall from [11] that an element $\Phi$ of
$\per{X \times Y}$ gives rise to an integral transform $\Phi$ from
$\per{X}$ to $\per{Y}$. This is a morphism in
$\text{Ho}(\text{dg-cat})$, the category of DG-categories modulo
quasi-equivalences. The functor from $\text{D}(\per{X})$ to
$\text{D}(\per{Y})$ induced by $\Phi$ is the functor $$E \mapsto
\pi_{Y*}(\pi_X^*E \otimes^{\mathbb L} \Phi)\text{ . }$$ $\Phi$
induces a map from $\hochh{\per{X}}$ to $\hochh{\per{Y}}$ and hence,
a map from $\hochh{X}$ to $\hochh{Y}$ which we shall denote by
$\natl{\Phi}$. We now state the following consequence of Theorem 4.
Like Theorem 4, Theorem 5 implies that $\natl{\Phi}$ depends only on
the image of $\Phi$ in $\text{D}(\per{X \times
Y})$.\\

\begin{thm} For any $\Phi$ in $\per{X \times Y}$,
$$\natl{\Phi}(x)=\langle x,\text{Ch}(\Phi)\rangle_{\text{Shk}} \in \hochh{Y} $$
for all $x \in \hochh{X}$.
\end{thm}

\textbf{Sketch of proof of Theorem 5.} Theorem 5 is a direct
consequence of Theorem 4 and the work of B. Toen [11]. Given two
DG-categories $\calg C$ and $\calg D$, [11] constructs a DG-category
$\text{RHom}(\calg C,\calg D)$. Let $X$ and $Y$ be quasi compact
separated schemes over $\mathbb K$. Let $E$ and $F$ be compact
generators of $\text{D}_{qcoh}(X)$ and $\text{D}_{qcoh}(Y)$. Recall
that in [11] it was shown that there is an identification
$$\beta:\per{{\calg A}(E)^{\text{op}} \otimes {\calg A}(F)}  \rar \text{RHom}(\per{{\calg A}(E)},\per{{\calg A}(F)})$$
$$\Phi \leadsto M \mapsto M \otimes_A \Phi $$
in $\text{Ho}(\text{dg-cat})$. Similarly, there is an identification
$$\gamma:\per{X \times Y} \rar \text{RHom}(\per{X},\per{Y})$$
in $\text{Ho}(\text{dg-cat})$. If $\Phi$ is in $\per{X \times Y}$,
$\gamma(\Phi)$ is the integral transform $\Phi$ from $\per{X}$ to
$\per{Y}$ that we described before stating Theorem 5. We abuse
notation and use $\eta$ to denote the quasi-equivalences
$\per{{\calg A}(E)} \simeq \per{X}$,$\per{{\calg A}(E)^{\text{op}}
\otimes {\calg A}(F)} \simeq \per{X \times Y}$  and\\
$\text{RHom}(\per{{\calg A}(E)},\per{\calg A(F)}) \simeq \text{RHom}(\per{X},\per{Y})$ described in [11].  \\

It was shown in Section 8 of [11] that the following diagram
commutes in $\text{Ho}(\text{dg-cat})$.
$$\begin{CD}
\per{X \times Y} @>\eta^{-1}>> \per{{\calg A}(E)^{\text{op}} \otimes {\calg A}(F)} \\
@VV{\gamma}V       @V{\beta}VV\\
\text{RHom}(\per{X},\per{Y}) @>\eta^{-1}>> \text{RHom}(\per{{\calg
A}(E)},\per{{\calg A}(F)})
\\
\end{CD}$$
Theorem 5 is then a direct consequence of Theorem 4 and the above commutative diagram.\\

\textbf{Remark.} Instead of choosing a compact generator $E$ of $\text{D}_{qcoh}(X)$ and using the DG-algebra ${\calg A}(E)$ to define $\langle \text{ },\text{ }\rangle_{\text{Shk}}$ on $\hochh{X}$, we could make do with any DG-algebra $A$ such that $\per{A}$ is quasi-equivalent to $\per{X}$. \\

\section{The Mukai pairing.}

\subsection{Some recollections.}
Let $X$ be a smooth proper scheme. Let $S_X$ denote the shifted line
bundle on $X$ tensoring with which yields the Serre duality functor
on the bounded derived category $\text{D}^{b}(X)$ of coherent
$\strc$-modules. If $f:X \rar Y$ is a morphism of schemes,
$f_*$,$f^*$ etc shall denote the corresponding derived functors in
this section. Let $\Delta:X \rar X \times X$ denote the diagonal
embedding. Let $\Delta_!$ denote the left adjoint of $\Delta^*$. Let
$\strcc{\Delta}$ denote $\Delta_* \strc$. Recall from [1] that there
is an isomorphism
$$\hochh{X} \simeq \text{RHom}_{X \times X}(\Delta_!\strc, \strcc{\Delta}) \text{ . }$$
Since $\Delta_!\strc \simeq \Delta_*S_X^{-1}$, tensoring with $\pi_2^*S_X$ yields an isomorphism
$$D:\rhh(\Delta_!\strc,\Delta_*\strc) \rar \rhh(\Delta_*\strc,\Delta_*S_X) \text{ . }$$
\textbf{Definition.}The Mukai pairing $\langle \text{ },\text{ } \rangle_M$ on $\hochh{X}$ is the pairing
$$v \otimes w \leadsto \text{tr}_{X \times X}(D(v) \circ w) $$
  where $\text{tr}_{X \times X}$ denotes the Serre duality trace on $X \times X$. The same pairing was constructed in the DG-algebra
  setup in [9]. \\

Recall that the Hochschild-Kostant-Rosenberg map $\hkr$ induces an isomorphism
$$\hoch{i}{X} \simeq \oplus_i \text{H}^{j-i}(X,\Omega_X^{j})$$ which we shall also denote by $\hkr$.
Let $\int_X$ denote the linear functional on $\oplus_{p,q} \text{H}^p(X,\Omega_X^q)$ that coincides with
the Serre duality trace on $\text{H}^n(X,\Omega_X^n)$ and vanishes on other direct summands.
 Let $*$ denote the involution on $\oplus_{p,q}\text{H}^{p}(X,\Omega_X^q)$ that acts on the summand
 $\text{H}^p(X,\Omega_X^q)$ by ${(-1)}^p$.
 The following result (implicitly in [7] and explicitly in [8]) computes $\gpair_M$ at the level of Hodge cohomology.

\begin{thm} For $a,b \in \hochh{X}$,
$$\spair{a}{b}_M = \int_X \hkr(a)^*\hkr(b)\text{td}(T_X) \text{ . }$$
\end{thm}

\subsection{Integral transforms in Hochschild homology.}

Any $\Phi \in \per{X \times Y}$ yields an integral transform
$$\Phi:\per{X} \rar \per{Y} $$
as described in Section 1.2. Note that if $\Psi \in \per{Y \times
Z}$,the image of the kernel of the integral transform $\Psi \circ
\Phi$ in $\text{D}(\per{X \times Z})$ is precisely
$\pi_{XZ*}(\pi_{XY}^* \Phi \otimes^{\mathbb L} \pi_{YZ}^*\Psi)$. A
priori, there is more than one construction of the corresponding
integral transform
$\Phi_*:\hochh{X} \rar \hochh{Y}$ such that \\
a. $(\Psi \circ \Phi )_*=\Psi_* \circ \Phi_*$. \\
b. The following diagram commutes.
$$\begin{CD}
\text{D}(\per{X}) @>\Phi>> \text{D}(\per{Y})\\
@VV{\text{Ch}}V               @V{\text{Ch}}VV\\
\hoch{0}{X} @>\Phi_*>> \hoch{0}{Y}\\
\end{CD} $$
For example, $\natl{\Phi}$ is seen to satisfy these properties
without much difficulty. Another construction of $\Phi_*$ was given
by A. Caldararu in [1]. Broadly speaking, one views $\hochh{X}$ as
an "ext of functors", $\text{Ext}(S_X^{-1},\text{id})$. This can be
done rigorously as in [3]. Let $\Phi^{\vee}$ be a left adjoint of
$\Phi$. Then, if $\alpha \in \text{Ext}(S_X^{-1},\text{id})$,
$\Phi_*(\alpha)$ is the following composite where the unlabeled
arrows are adjunctions.
$$\begin{CD} S_Y^{-1} \\
@VVV\\
\Phi \circ \Phi^{\vee} \circ S_Y^{-1}\\
 @V{=}VV\\
 \Phi \circ S_X^{-1} \circ S_X \circ \Phi^{\vee} \circ S_Y^{-1} @>{\text{id} \circ \alpha \circ \text{id} \circ \text{id}\circ \text{id}}>> \Phi \circ S_X \circ \Phi^{\vee} \circ S_Y^{-1} @>>> \text{id}_Y \\\end{CD} $$
Theorem 6 enables us to give yet another construction of $\Phi_*$.
This construction of $\Phi_*$ is motivated by Theorem 4,and plays a
key role in relating the Mukai pairing to the natural pairing
constructed in Section 1. In the rest of this section, the
identification of $\hochh{X \times Y}$ with $\hochh{X} \otimes
\hochh{Y}$ will be via the inverse of the Kunneth isomorphism.
Recall that if $\Phi \in \per{X \times Y}$, the Chern character
$\text{Ch}(\Phi) \in \hochh{X \times Y} \simeq \hochh{X} \otimes
\hochh{Y}$ may be viewed as a $\mathbb K$-linear map from $\mathbb
K$ to $\hochh{X} \otimes \hochh{Y}$.
 Let $\text{W}:\hochh{X} \rar \hochh{X}$ be the unique involution corresponding via $\hkr$ to the involution $*$ on
  $\oplus_{p,q} \text{H}^p(X,\Omega_X^q)$ mentioned in the previous subsection. \\

$\textbf{Construction.}$ We define $\muk{\Phi}:\hochh{X} \rar \hochh{Y}$ to be the composite
$$\begin{CD} \hochh{X} \\\
@V {\text{id}\otimes \text{Ch}(\Phi) }VV\\
   \hochh{X}^{\otimes 2} \otimes \hochh{Y}  @> \text{W} \otimes \text{id} \otimes \text{id}>>  \hochh{X}^{\otimes 2} \otimes \hochh{Y}  @> \gpair_M \otimes \text{id}>> \hochh{Y}\\ \end{CD} $$

\begin{prop}
If $\Phi \in \per{X \times Y}$ and $\Psi \in \per{Y \times Z}$ then
$$\muk{(\Psi \circ \Phi)} = \muk{\Psi} \circ \muk{\Phi} \text{ . }$$
\end{prop}

\begin{proof} We shall denote $\oplus{p,q}\text{H}^p(X,\Omega_X^q)$ by $\text{H}^{\bullet}(X)$.
Recall ( Theorem 4.5 \\ in [2]) that for any smooth scheme $Z$,
$\hkr \circ \text{Ch} = \text{ch}$,
the right hand side being the familiar Chern character map from $\text{D}(\per{Z})$ to $\text{H}^{\bullet}(Z)$. \\

Let $a \in \hochh{X}$. Note that $\hoch{0}{X \times Y} \simeq \oplus_i \hoch{i}{X} \otimes \hoch{-i}{Y}$. Hence,
$$\text{Ch}(\Phi) = \sum_i \sum_{\lambda(i) \in I_i} \alpha_{\lambda(i)} \otimes \beta_{\lambda(i)} $$ for some index sets $I_i$ and $\alpha_{\lambda(i)} \in \hoch{i}{X}$ and $\beta_{\lambda(i)} \in \hoch{-i}{Y}$.
By Theorem 6 and the construction of $\muk{\Phi}$,
\begin{equation} \label{mukpush} \hkr(\muk{\Phi}(a)) = \sum_i \sum_{\lambda(i) \in I_i} (\int_X  \hkr(a) \hkr(\alpha_{\lambda(i)})\text{td}(T_X))  \hkr(\beta_{\lambda(i)}) \text{ . }\end{equation}

Now suppose that $$\text{Ch}(\Psi) = \sum_j \sum_{\mu(j) \in J_j}
\gamma_{\mu(j)} \otimes \delta_{\mu(j)} $$ for some index sets $J_j$
and $\gamma_{\mu(j)} \in \hoch{j}{Y}$ and $\delta_{\mu(j)} \in
\hoch{-j}{Z}$. Then, by $\eqref{mukpush}$,
$$\hkr(\muk{\Psi} \circ \muk{\Phi}(a)) = $$ $$ \sum_{i,j}
\sum_{\lambda(i) \in I_i,\mu(j) \in J_j} \hkr(\delta_{\mu(j)})
(\int_X \hkr(a) \hkr(\alpha_{\lambda(i)})  \text{td}(T_X)) (\int_Y
\hkr(\beta_{\lambda(i)}) \hkr(\gamma_{\mu(j)})  \text{td}(T_Y))
\text{ . }$$

Recall that $\Psi \circ \Phi = \pi_{XZ*}(\pi_{YZ}^*\Psi \otimes
\pi_{XY}^*\Phi)$ The desired proposition will follow from
$\eqref{mukpush}$ if we can show that \begin{equation}
\label{toshow} \text{ch}(\Psi \circ \Phi) = \sum_{i,j}
\sum_{\lambda(i) \in I_i,\mu(j) \in J_j} (\int_Y
\hkr(\beta_{\lambda(i)}) \hkr(\gamma_{\mu(j)}) \text{td}(T_Y))
\alpha_{\lambda(i)} \otimes \delta_{\mu(j)}
 \text{ . }\end{equation}

Recall that after identifying $\text{H}^{\bullet}(X \times Y)$ with
$\text{H}^{\bullet}(X) \otimes \text{H}^{\bullet}(Y)$, $\pi_{Y*}$
gets identified with $\int_X \otimes \text{id}$. Also, $\pi_Y^*$ is
identified with the map $a \leadsto 1 \otimes a$ from
$\text{H}^{\bullet}(Y)$ to $\text{H}^{\bullet}(X \times Y)$. With
this in mind, $\eqref{toshow}$ can be rewritten as,
$$ \text{ch}(\Psi \circ \Phi) = \pi_{XZ*} (\text{ch}(\pi_{XY}^*(\Phi)).\text{ch}(\pi_{YZ}^*
\Psi).\pi_{Y}^*\text{td}(T_Y)) \text{ . }$$ This follows directly
from the Riemann-Roch-Hirzebruch theorem
applied to the map $\pi_{XZ}:X \times Y \times Z \rar X \times Z$.\\

\end{proof}

Let $\strcc{\Delta} = \Delta_* \strc$ be treated as the kernel of an
integral transform from $X$ to $X$. Then,

\begin{prop}
$$\muk{\strcc{\Delta}} = \text{id} \text{ . }$$
\end{prop}

\begin{proof}

Since $\strcc{\Delta}$ is the kernel of the identity,
$\strcc{\Delta} \circ \strcc{\Delta} = \strcc{\Delta}$. By
Proposition 1, $\muk{\strcc{\Delta}}$ is an idempotent endomorphism
of $\hochh{X}$. To prove that it is the identity, it suffices to
show that it is surjective.\\

For this, note that $\natl{\strcc{\Delta}} = \text{id}$. By theorem
4,
$$\text{Ch}(\strcc{\Delta}) = \sum_i \sum_k e_{i,k} \otimes
f_{i,k}$$ where the $e_{i,k}$ form a basis of $\hoch{i}{X}$ and the
$f_{i,k}$ form a basis of $\hoch{-i}{X}$ such that
$$\spair{f_{i,k}}{e_{i,l}}_{\text{Shk}} = \delta_{k,l} \text{ . }$$
The $\delta$ on the right hand side of the above equation is the
Kronecker delta.\\

Let $\text{W}$ be the involution on $\hochh{X}$ which we defined
earlier before constructing $\muk{\Phi}$. It follows that if $x \in
\hoch{i}{X}$, then
$$\muk{\strcc{\Delta}}(x) = \sum_k \spair{\text{W}(x)}{e_{-i,k}}_M f_{-i,k} \text{
. }$$ Recall from [1] that the pairing $\gpair_M$ is non-degenerate.
Moreover, $\text{W}$ is an involution on $\hochh{X}$. Since the
$e_{-i,k}$ form a basis of $\hoch{-i}{X}$, there exist elements
$x_k$ in $\hoch{i}{X}$ such that
$$\spair{\text{W}(x_l)}{e_{-i,k}}_M = \delta_{kl} \text{ . }$$
Clearly, $$\muk{\strcc{\Delta}}(x_k) = f_{-i,k} \text{ . }$$ This
proves that $\muk{\strcc{\Delta}}$ is surjective, as was desired.
\end{proof}

We are now ready to prove Theorem 1.

\subsection*{Proof of Theorem 1.}
\begin{proof}
This follows almost immediately from the fact that $\muk{\strcc{
\Delta}} = \natl{\strcc{\Delta}} = \text{id}: \hochh{X} \rar
\hochh{X}$. Since $\natl{\strcc{\Delta}} = \text{id}$,
$$\spair{f_{-i,k}}{e_{-i,l}}_{\text{Shk}}=\delta_{k,l} \text{ . }$$
On the other hand since $\muk{\strcc{\Delta}} = \text{id}$ by
Proposition 2,
$$\spair{\text{W}(f_{-i,k})}{e_{-i,l}}_M =\delta_{k,l} \text{ . }$$
It follows from the $\mathbb K$ bi-linearity of the pairings
$$(a,b) \leadsto \spair{a}{b}_{\text{Shk}} $$
$$(a,b) \leadsto \spair{\text{W}(a)}{b}_M $$ that
\begin{equation} \label{mainresult}
\spair{a}{b}_{\text{Shk}}=\spair{\text{W}(a)}{b}_M \text{ .}
\end{equation}
Recall that $\vee$ denotes the involution on $\hochh{X}$
corresponding via $\hkr$ to the involution on
$\text{H}^{\bullet}(X)$ which acts on the direct summand
$\text{H}^{q}(X,\Omega_X^p)$ by multiplication by ${(-1)}^p$. Now,
$\hkr(a)^* \cup \hkr(b) = \hkr(b) \cup \hkr(a^{\vee})$ in
$\text{H}^{\bullet}(X)$. Hence, Theorem 6 may be rewritten to say
that
$$\spair{a}{b}_M = \int_X \hkr(b)\hkr(a^{\vee}) \text{td}(T_X) \text{ . } $$
By $\eqref{mainresult}$,
$$ \spair{a}{b}_{\text{Shk}} = \int_X \hkr(a)\hkr(b)\text{td}(T_X) =
\spair{b^{\vee}}{a}_M \text{ . }$$

This proves Theorem 1.\\
\end{proof}

Recall from [1] that the integral transform from $\text{D}(\per{Y})$
to $\text{D}(\per{X})$ due to $\rhhh(\Phi,\strcc{X \times Y})
\otimes^{\mathbb L} \pi_X^*S_X$ is the right adjoint of that from
$\text{D}(\per{X})$ to $\text{D}(\per{Y})$ due to $\Phi$. Let
$\Phi^!$ denote $\rhhh(\Phi,\strcc{X \times Y}) \otimes^{\mathbb L}
\pi_X^*S_X$. We also have the following proposition, which shows
that $\muk{\Phi}$ is a "good candidate" for the integral transform
on Hochschild homology
defined by $\Phi$.\\

\begin{prop}
If $x \in \hochh{X}$ and $y \in \hochh{Y}$, then
$$\spair{\muk{\Phi}(x)}{y}_M = \spair{x}{\muk{\Phi^!}(y)}_M \text{ .
}$$
\end{prop}

\begin{proof} The notation used in this proof is as in the proof of
Proposition 1.  Assume that after identifying $\hochh{X \times Y}$
with $\hochh{X} \otimes \hochh{Y}$ (via the inverse of the Kunneth
map),
$$\text{Ch}(\Phi) = \sum_i \sum_{\lambda(i) \in I_i} \alpha_{\lambda(i)} \otimes \beta_{\lambda(i)} $$ for some
 index sets $I_i$ and $\alpha_{\lambda(i)} \in \hoch{i}{X}$ and $\beta_{\lambda(i)} \in \hoch{-i}{Y}$.
 Then, by Theorem 6 and $\eqref{mukpush}$,
$$\spair{\muk{\Phi}(x)}{y}_M = $$ $$ \sum_i \sum_{\lambda(i) \in I_i}
(\int_X \hkr(x) \hkr( \alpha_{\lambda(i)})\text{td}(T_X))(\int_Y
\hkr(\beta_{\lambda(i)})^* \hkr(y) \text{td}(T_Y)) \text{ . }$$

Note that $\text{Ch}(\Phi^!) = \sum_i \sum_{\lambda(i) \in I_i}
{(-1)}^i \text{W}(\beta_{\lambda(i)}) \otimes
[\text{W}(\alpha_{\lambda(i)}).\text{Ch}(S_X)]$. The ${(-1)}^i$
comes from the fact that the composite
$$\begin{CD} \hochh{X} \otimes \hochh{Y} @>K>> \hochh{X \times Y} @>{K}^{-1}>> \hochh{Y} \otimes \hochh{X} \end{CD} $$
is the signed map swapping factors. It follows from Theorem 6 and
$\eqref{mukpush}$ that
$$\spair{x}{\muk{\Phi^!}(y)}_M = $$ $$ \sum_i \sum_{\lambda(i) \in I_i} {(-1)}^i (\int_Y
\hkr(\hkr(y) \beta_{\lambda(i)})^*  \text{td}(T_Y))(\int_X \hkr(x)^*
\hkr(\alpha_{\lambda(i)})^* \text{ch}(S_X) \text{td}(T_X)) $$
$$= \sum_i \sum_{\lambda(i) \in I_i}  (\int_Y
  \hkr(\beta_{\lambda(i)})^* \hkr(y)  \text{td}(T_Y))(\int_X \hkr(x)^*
\hkr(\alpha_{\lambda(i)})^* \text{ch}(S_X) \text{td}(T_X)) \text{ .
}$$

Now, if $n$ is the dimension of $X$, $\text{ch}(S_X) =
{(-1)}^n\text{ch}(\Omega_X^n)$. Also,
$\text{td}(T_X)\text{ch}(\Omega_X^n) = \text{td}(T_X)^*$ (see [2]).
It follows that
$$ \hkr(x)^*
\hkr(\alpha_{\lambda(i)})^* \text{ch}(S_X) \text{td}(T_X)) =
{(-1)}^n (\hkr(x) \hkr( \alpha_{\lambda(i)})\text{td}(T_X))^* \text{
. }$$ Hence,$$\int_X \hkr(x)^* \hkr(\alpha_{\lambda(i)})^*
\text{ch}(S_X) \text{td}(T_X) = \int_X \hkr(x) \hkr(
\alpha_{\lambda(i)})\text{td}(T_X) $$This proves the desired
proposition.

\end{proof}

Note that Proposition 1 and Proposition 3 parallel Theorems 5.3 and
7.3 respectively in [1]. However, since we use the Riemann-Roch
theorem for (proper) projections to prove Proposition 1, the
construction of $\muk{\Phi}$ by itself does not amount to a
self-contained construction of integral transforms in Hochschild
homology at this stage. However, it helps prove Theorem 1, which in
turn leads to Theorem 2, showing that all three constructions of
integral transforms in Hochschild homology coincide. In particular,
it tells us that the integral transform constructed by A. Caldararu
[1] coincides with the more "natural" construction
of the integral transform constructed by D. Shklyarov [10].\\

Let $\Phi \in \per{X \times Y}$. Denote the integral transform
$\Phi_*:\hochh{X} \rar \hochh{Y}$ constructed by A. Caldararu [1]
and described briefly earlier in this section by $\cald{\Phi}$.\\

\subsection*{Proof of Theorem 2.}

\begin{proof}
That $\muk{\Phi}=\natl{\Phi}$ is an immediate consequence of Theorem
1 and Theorem 5. We therefore need to show that $\muk{\Phi} =
\cald{\Phi}$. For this, we will follow D. Shklyarov and imitate the
proof of Theorem 4 (Theorem 3.4 in [10]) in [10]).\\

{\it Step 1:}Recall that if $\Phi \in \per{X \times Y}$ and $\Phi'
\in \per{X' \times Y'}$ , $\Phi \boxtimes \Phi' \in \per{X \times
X'\times Y \times Y'}$. We then have integral transforms in
Hochschild homology
$$\muk{\Phi}:\hochh{X} \rar \hochh{Y} , \muk{\Phi'}:\hochh{X'}
\rar \hochh{Y'} $$ $$ \muk{(\Phi \boxtimes \Phi')}:\hochh{X \times
X'} \rar \hochh{Y \times Y'} \text{ . }$$ Identify $\hochh{X \times
X'}$ and $\hochh{Y \times Y'}$ with $\hochh{X} \otimes \hochh{X'}$
and $\hochh{Y} \otimes \hochh{Y'}$ respectively via the inverse of
the relevant Kunneth isomorphisms. It follows from the construction
of $\muk{\Phi}$ that
$$\muk{(\Phi \boxtimes \Phi')} = \muk{\Phi} \otimes \muk{\Phi'} \text{
. }$$

Similarly, we have integral transforms in Hochschild homology
$$\cald{\Phi}:\hochh{X} \rar \hochh{Y} , \cald{\Phi'}:\hochh{X'}
\rar \hochh{Y'} $$ $$ \cald{(\Phi \boxtimes \Phi')}:\hochh{X \times
X'} \rar \hochh{Y \times Y'} \text{ . } $$ It can be verified
without much difficulty (see [16], Lemma 2.1 for instance) that
$$\cald{(\Phi \boxtimes \Phi')} = \cald{\Phi} \otimes \cald{\Phi'} \text{
. }$$

{\it Step 2:} Note that $\Phi \in \per{X \times Y}$ may also be
thought of as the kernel of an integral transform from $\text{ Spec
} \mathbb K$ to $X \times Y$. We will denote $\Phi$ thought of in
this manner by $\Phi_{\text{pt} \rar X \times Y}$. Let $\Delta$
denote $\strcc{\Delta}$ thought of as the kernel of an integral
transform from $X \times X$ to $\text{Spec }\mathbb K$. Also
identify $\hochh{X}$ with $\hochh{X} \otimes \hochh{\text{Spec }
\mathbb K}$ via the map $y \leadsto y \otimes 1$. Then,
$$\Phi = \Delta \circ (\strcc{\Delta} \boxtimes \splt) $$
$$\implies \muk{\Phi} = \muk{\Delta} \circ \muk{(\strcc{\Delta} \boxtimes
\splt)} = \muk{\Delta} \circ (\muk{\strcc{\Delta}} \otimes
\muk{(\splt)}(1)) $$ $$ \cald{\Phi} = \cald{\Delta} \circ
\cald{(\strcc{\Delta} \boxtimes \splt)} = \cald{\Delta} \circ
(\cald{\strcc{\Delta}} \otimes \cald{(\splt)}(1)) $$ Now, by
Proposition 2,
$$\cald{\strcc{\Delta}}=\muk{\strcc{\Delta}} = \text{id} \text{ .
}$$ Also, $\cald{(\splt)}(1) = \text{Ch}(\Phi)$ by Definition 6.1 in
[1] and Theorem 4.5 in [2]. $\muk{(\splt)}(1) = \text{Ch}(\Phi)$ by
the construction of $\muk{(\splt)}$. We therefore , need to show
that
$$\muk{\Delta}=\cald{\Delta}: \hochh{X \times X} \rar \hochh{\text{Spec }\mathbb K} = \mathbb K \text{ . }$$
With the above identification of $\hochh{\text{Spec }\mathbb K}$
with $\mathbb K$, for any $x \in \hochh{\text{Spec
}\mathbb K}$, $x = \spair{x}{1}_M$. Let $\Delta^!$ denote \\
$\rhhh(\Delta, \strcc{X \times X}) \otimes^{\mathbb L} S_{X \times
X}$.  If $\alpha \in \hochh{X \times X}$,
$$\spair{\muk{\Delta}(\alpha)}{1}_M =
\spair{\alpha}{\muk{\Delta^!}(1)}_M $$ by Proposition 3. By Theorem
7.3 in [1], $$\spair{\cald{\Delta}(\alpha)}{1}_M =
\spair{\alpha}{\cald{\Delta^!}(1)}_M \text{ . }$$ Now,
$\cald{\Delta^!}(1)= \text{Ch}(\Delta^!)$ by Definition 6.1 in [1]
and Theorem 4.5 in [2]. $\muk{\Delta^!}(1)= \text{Ch}(\Delta^!)$ by
the construction of $\muk{\Delta^!}$. This yields the
desired theorem.\\
\end{proof}

\subsection{When $X$ is Calabi-Yau.}

In such a situation, $\text{D}^{b}(X)$ can be be thought of as the
category of open states of the B-Model on $X$ (see [3]). The
corresponding algebra of closed states is the Hochschild cohomology
$\text{HH}^{\bullet}(\per{X}) \simeq \text{HH}^{\bullet}(X)$. As $X$
is Calabi-Yau, there is an identification
$$\text{HH}^{\bullet}(X) \simeq \hochh{X} \text{ .}$$ The Mukai pairing constructed by A. Caldararu in [1] on $\hochh{X}$
then gives a pairing on $\text{HH}^{\bullet}(X)$. Moreover, for any
${\mathcal E} \in \text{D}^{b}(X)$, there are natural maps
$$\iota^{\calg E} : \text{Hom}_{\dcat{X}}(\calg E,\calg E) \rar \text{HH}^{\bullet}(X) $$
$$\iota_{\mathcal E}:\text{HH}^{\bullet}(X) \rar \text{Hom}_{\dcat{X}}(\calg E,\calg E) $$
as constructed in [3]. The Cardy condition verifies that this data
gives a topological quantum field theory. Of course, the Mukai
pairing in this case is the pairing obtained by the action of the
class of a genus $0$ Riemann-surface with two incoming closed
boundaries and no outgoing boundary in $\text{H}_{0}({\mathcal
M}_0(2,0))$ on $\hochh{X}$, the action coming from the fact that
$\text{HH}^{\bullet}(X)$ with Mukai pairing
is a "good" algebra of closed states as verified by the Cardy condition.\\

On the other hand, [4]  gives the category of open states of the
B-Model on $X$ as an $A_{\infty}$ enrichment of $\dcat{X}$. The
closed TCFT one associates with this category has homology
$$\hochh{X} \simeq \text{HH}^{\bullet}(X) \text{ . }$$ This is also
equipped with a pairing coming out of the action of the class of a
genus $0$ Riemann-surface with two incoming closed boundaries and no
outgoing boundary in $\text{H}_{0}({\mathcal M}_0(2,0))$ on the
homology of the closed TCFT one constructs in [4] from the B-Model.
 Whether these pairings coincide
is however, not clear currently.\\

Theorem 1 is similar Conjecture 6.2 in [10] for Calabi-Yau algebras
$A$ such that $\per{A}$ is quasi-equivalent to $\per{X}$ for some
quasi-compact separated smooth scheme $X$.\\

\subsection{Proof of Theorem 3.}

\textbf{The sheafification of the Dennis trace map.} Let us briefly
recall how the sheafification of the Dennis trace map is
constructed. The material we are recalling is from [12],[13],[14]
and [15]. Let $X$ be a smooth quasicompact separated scheme. As in
Section 1.2, choose a compact generator $E$ of $\text{D}_{qcoh}(X)$
and a DG-algebra ${\calg A}(E)$ such that $\per{{\calg A}(E)}$ is
quasiequivalent to $\per{X}$. Let $\perr{{\calg A}(E)}$ be the exact
category whose objects are those of $\per{{\calg A}(E)}$ such that
$$\text{Hom}_{\perr{{\calg A}(E)}}(M,N) = \text{Z}^{0}(\text{Hom}_{\per{{\calg
A}(E)}}(M,N)) \text{ . }$$ As pointed out by B. Keller in [14],
using the Waldhausen structure of $\perr{{\calg A}(E)}$, we can
construct a Dennis trace map
$$\text{Dtr}: K_i(X) \simeq K_i(\perr{{\calg A}(E)}) \rar
\text{HH}_{i,\text{McC}}(\perr{{\calg A}(E)}) \text{  } \forall i
\geq 0 \text{ . }$$ Here, $\text{HH}_{i,\text{McC}}$ is the
Hochschild homology constructed by R. McCarthy in [15]. As Keller
further points out in [14], there is a natural transformation
$$\text{HH}_{i,\text{McC}}(\perr{{\calg A}(E)}) \rar \hoch{i}{\perr{{\calg
A}(E)}} \text{ . }$$ Further, we also have a natural transformation
$$\hoch{i}{\perr{{\calg A}(E)}} \rar \hoch{i}{\per{{\calg A}(E)}}
\text{ . }$$ The obvious compositions then give us a map
$$\text{Ch}^i: K_i(X) \simeq K_i(\perr{{\calg A}(E)}) \rar \hoch{i}{\per{{\calg
A}(E)}} \simeq \hoch{i}{X} \text{ . }$$

Let $Y$ be a smooth quasicompact separated scheme. Let $F$ and
${\calg A}(F)$ be as in Section 1.2. Let $\Psi \in \per{{\calg
A}(F)^{\text{op}} \otimes {\calg A}(E)}$. The following proposition,
analogous to Theorem 7.1 of [1], says that the sheafification of the
Dennis trace map is "functorial".

\begin{prop}
The following diagram commutes.
$$\begin{CD}
K_i(\perr{{\calg A}(E)}) @>\Psi_*>> K_i(\perr{{\calg A}(F)}) \\
@VV{\text{Ch}^i}V      @V{\text{Ch}^i}VV\\
\hoch{i}{\per{{\calg A}(E)}} @>\natl{\Psi}>> \hoch{i}{\per{{\calg
A}(F)}}\\
\end{CD}$$

\end{prop}

\begin{proof} This proposition will follow easily once we verify
that \\ $\Psi : \perr{{\calg A}(E)} \rar \perr{{\calg A}(F)}$
preserves cofibrations and weak equivalences. By [12], the weak
equivalences in $\perr{A}$ for any DG-algebra $A$ are
quasiisomorphisms. The cofibrations in $\perr{A}$ are morphisms of
$A$-modules that admit retractions as morphisms of graded
$A$-modules. That $\Psi$ preserves cofibrations follows without
difficulty from the fact that $\Psi:\per{{\calg A}(E)} \rar
\per{{\calg A}(F)}$ is a DG-functor. That $\Psi$ preserves weak
equivalences follows from the fact that perfect modules are
homotopically projective (see Proposition 2.5 of [10]).

\end{proof}

\textbf{Proof of Theorem 3.} We warn the reader that in the proof
that follows, $X$ and $Y$ denote {\it proper} smooth quasicompact
separated schemes.
\begin{proof}

{\it Step 1:} Let $\Phi \in \per{X \times Y}$.The first step is to
note that even though $Z$ is not necessarily proper, the kernel
$\Phi \boxtimes \strcc{\Delta_Z} \in \per{X \times Z \times Y \times
Z}$ induces an integral transform from $\per{X \times Z}$ to $\per{Y
\times Z}$. This follows from the fact that if $E$ and $F$ are
compact generators of $\text{D}_{\text{coh}}(X)$ and
$\text{D}_{\text{coh}}(Z)$ respectively, the compact generator $E
\boxtimes F :=\pi_X^*E \otimes \pi_Z^*F$ of $\text{D}_{\text{coh}}(X
\times Z)$ is mapped by the integral transform with kernel $ \Phi
\boxtimes \strcc{\Delta_Z} $ to the perfect complex $ \pi_{Y*}(\Phi
\otimes^{\mathbb L} \pi_X^*E) \boxtimes F
$.\\

Also, after identifying $\hochh{X \times Z}$ and $\hochh{Y \times
Z}$ with $\hochh{X} \otimes \hochh{Z}$ and $\hochh{Y} \otimes
\hochh{Z}$ respectively via the inverse of the relevant Kunneth
isomorphisms,
\begin{equation} \label{kunn} \natl{(\Phi \boxtimes \strcc{\Delta_Z})} = \natl{\Phi} \otimes
\text{id} :\hochh{X} \otimes \hochh{Z} \rar \hochh{Y} \otimes
\hochh{Z} \text{ . }\end{equation} This follows from the facts that
$\natl{\strcc{\Delta_Z}} = \text{id}$ and from Proposition 2.11 of
[10].\\

{\it Step 2:} By the Proposition 4, the following diagram commutes.
\begin{equation} \label{step1}
\begin{CD}
K_i(\per{X \times Z}) @>{(\Phi \boxtimes \strcc{\Delta_Z} )_*}>> K_i(\per{Y \times Z})\\
@VV{\text{Ch}^i}V       @V{\text{Ch}^i}VV\\
\hoch{i}{X \times Z} @>\natl{(\Phi \boxtimes \strcc{\Delta_Z})}>> \hoch{i}{Y \times Z} \\
\end{CD} \end{equation}

After identifying $\hochh{X \times Z}$ and $\hochh{Y \times Z}$ with
$\hochh{X} \otimes \hochh{Z}$ and $\hochh{Y} \otimes \hochh{Z}$
respectively via the inverse of the relevant Kunneth isomorphisms,we
have the following commutative diagram by $\eqref{step1}$ and
$\eqref{kunn}$.

\begin{equation} \label{step5}
\begin{CD}
K_i(\per{X \times Z}) @>{(\Phi \boxtimes \strcc{\Delta_Z})_*}>> K_i(\per{Y \times Z})\\
@VV{\text{Ch}^i}V       @V{\text{Ch}^i}VV\\
\oplus_{p+q=i} \hoch{p}{X} \otimes \hoch{q}{Z} @>\natl{\Phi} \otimes \text{id} >> \oplus_{p+q=i}\hoch{p}{Y} \otimes \hoch{q}{Z} \\
\end{CD} \end{equation}

Now, it follows from Theorem 1 and Theorem 3 that
$\muk{\Phi}=\natl{\Phi}$. Hence, by $\eqref{step5}$ and Proposition
3, \begin{equation} \label{step21} (\gpair_M \otimes
\text{id}_{\hochh{Z}})( f^*y \otimes \text{Ch}^i(\alpha)) =
(\gpair_M \otimes \text{id}_{\hochh{Z}}) (y \otimes (\text{id}
\times f)_* \text{Ch}^i(\alpha) )
\end{equation} for any $\alpha \in K_i(Z \times X)$, $y \in
\hochh{Y}$. By Theorem 4, $\eqref{step21}$ can be rewritten to say
that
$$\int_X \hkr(f^*(y))^* \text{ch}^i(\alpha) \text{td}(T_X) = \int_Y \hkr(y)^* \text{ch}^i((f \times \text{id})_* \alpha) \text{td}(T_Y) $$
as elements of $\text{H}^{\bullet}(Z)$. The desired theorem now
follows from the facts that $f^*$ commutes with $\hkr$ (see Theorem
7 of [7]) and commutes with the involution $*$.

\end{proof}

\textbf{Address:} \\
Department of Mathematics\\
University of Oklahoma\\
Norman, OK-73019\\
\textbf{email:} aramadoss@math.ou.edu\\

\end{document}